\documentclass[11pt]{article}
\usepackage{xypic}

\newcommand{\eop}
{
\mbox{}
\nolinebreak
\hfill
\raisebox{-0.2mm}{\rule{2mm}{2mm}}
\newline
\newline
}

\newenvironment{lelemme}{\begin{flushleft} {\bf Lemme} \it} {\end{flushleft}}
\newenvironment{letheoreme}{\begin{flushleft}
                            {\bf Th{\'e}or{\`e}me} \it} {\end{flushleft}}

\newtheorem{theoreme}{Th{\'e}or{\`e}me}
\newtheorem{lemme}{Lemme}
\newtheorem{rem}{Remarque}

\newcommand{\Aa}{{\cal A}}
\newcommand{\Bb}{{\cal B}}
\newcommand{\Cc}{{\cal C}}
\newcommand{\Ee}{{\cal E}}
\newcommand{\Ff}{{\cal F}}
\newcommand{\Gg}{{\cal G}}
\newcommand{\Hh}{{\cal H}}
\newcommand{\Kk}{{\cal K}}
\newcommand{\Ii}{{\cal I}}
\newcommand{\pp}{{\bf P}}
\newcommand{\Oo}{{\cal O}}
\newcommand{\op}[2]{#1\Oo_{\pp^{#2}}}
\newcommand{\tp}[1]{\mbox{$T$}_{\pp^{#1}}}

\newcommand{\kk}{{\bf k}~}

\newcounter{para}
\newcommand{\paragraphe}[1]
{
  \stepcounter{para}
  \subsection*{
    \begin{center}
      \thepara.~#1
    \end{center}
  }
}

\title{Rang maximal pour $\tp{n}$}
\author{Fran{\c{c}}ois Lauze \\
        U.N.S.A Laboratoire de Math{\'e}matiques  J. Dieudonn{\'e},
        U.R.A 168\\
        Parc Valrose, 06108 Nice Cedex 2\\
        e-mail : lauze@math.unice.fr
       }
\date{}

\newcommand{\itemb}{\item[$\bullet$]}
\newcommand{\demo} {{\bf D{\'e}monstration : }}

\begin{document}
\maketitle
\begin{abstract}
  Let $\kk$ an algebraically closed field, $\pp^n$ the n-dimensional projective
space over $\kk$ and $\tp{n}$ the tangent vector bundle of $\pp^n$. In this
paper I prove the following result : for every integer $\ell$, for every
non-negative integer $s$, if $Z_s$ is the union of $s$ points in
sufficiently general position in $\pp^n$, then the restriction map
$H^0(\pp^n,\tp{n}(\ell))\rightarrow H^0(Z_s,\tp{n}(\ell)_{|Z_s})$ has
maximal rank. This result implies that the last non-trivial term of the
minimal free resolution of the homogeneous ideal of $Z_s$ is the
conjectured one by the Minimal Resolution Conjecture of Anna Lorenzini
(cf~\cite{mrc}).

\end{abstract}

\paragraphe{Introduction}
Soit {\bf k} un corps alg{\'e}briquement clos et $\pp^n = \pp^n_{\bf k}$
l'espace projectif de dimension $n$ sur {\bf k} et $S = {\bf k}[x_0,
\ldots,x_n]$ l'anneau des coordonn{\'e}es homog{\`e}nes sur $\pp^n$. Soit $a$
un entier non n{\'e}gatif et  $R := \{P_1,\ldots,P_a\} \subset \pp^n$ des
points en position suffisamment g{\'e}n{\'e}rale, $\Ii_R$ le faisceau
d'id{\'e}aux de $R$ et $I = \Gamma_*(\Ii_R)$ l'id{\'e}al homog{\`e}ne de $Z$
dans $S$. Si on pose $d := min\{\ell,h^0(\pp^n,\Ii_R(\ell)) > 0\}$, alors $I$
est engendr{\'e} par $I_d\oplus I_{d+1}$ o{\`u} $I_s$ d{\'e}signe la partie
homog{\`e}ne de degr{\'e} $s$ de $I$ (Th{\'e}orie de Castelnuovo-Mumford,
\cite{mumford}, p 99). La r{\'e}solution minimale libre de $I$ s'{\'e}crit
alors :
$$
  0\rightarrow L_n \rightarrow \ldots \rightarrow L_p \rightarrow \ldots
  \rightarrow L_0\rightarrow I\rightarrow 0
$$
avec
$$
   L_p = S(-d-p-1)^{a_p}\oplus S(-d-p)^{b_p}
$$
et
$$
   a_p := h^1(\pp^n,\Omega_{\pp^n}^{p+1}(d+p+1)\otimes \Ii_R)\;,\;
   b_p := h^0(\pp^n,\Omega_{\pp^n}^{p}(d+p)\otimes \Ii_R).
$$
Puisque l'application de restriction
\begin{equation}
  \label{eq1}
  H^0(\pp^n,\op{}{n}(\ell)) \rightarrow H^0(R,\Oo_R(\ell))
\end{equation}
est de rang maximal, c'est {\`a} dire surjective ou injective, les nombres
$b_0$,
$a_{n-1}$ et $b_n$ sont connus.

Les nombres $a_p$ et $b_p$ ont {\'e}t{\'e}s calcul{\'e}s dans~\cite{hirsim}
pour tout nombre $a$ de points suffisamment grand et tout entier p.
L'entier $a_n$ est {\'e}videmment nul et puisque l'application~(\ref{eq1})
est de rang maximal, on a
\begin{itemize}
\itemb $b_n = max\{0,h^0(\pp^n,\op{}{n}(d-1)) - a\} = 0$,
        par d{\'e}finition de $d$,
\itemb $a_{n-1} = max\{0,a-h^0(\pp^n,\op{}{n}(d-1))\} = a -
        \left(\!\!\!\begin{array}{c} n\!\!+\!\!d\!\!-\!\!1\\ n
        \end{array}\!\!\!\right)$, pour les m{\^e}mes raisons.
\end{itemize}
Le but de cet article est de calculer, sans restriction sur le nombre de
points, les entiers $a_{n-2}$ et $b_{n-1}$. On va {\'e}tablir le
th{\'e}or{\`e}me
suivant :\vspace{3mm}

\begin{theoreme}
  \label{th1}
  Avec les notations pr{\'e}c{\'e}dentes :
  \begin{itemize}
  \itemb $a_{n-2} = max\{0,na  - h^0(\pp^n,\Omega_{\pp^n}^{n-1}(d+n-1))\}$
  \itemb $b_{n-1} = max\{0,h^0(\pp^n,\Omega_{\pp^n}^{n-1}(d+n-1))-na\}$,
  \end{itemize}
  c'est {\`a} dire que l'application de restriction
  $$
    H^0(\pp^n,\Omega_{\pp^n}^{n-1}(d+n-1)) \rightarrow
    H^0(R,\Omega_{\pp^n}^{n-1}(d+n-1)|_R)
  $$
  est de rang maximal, i.e. surjective ou injective (puisque
  $rg(\Omega^{n-1}_{\pp^n})=n)$.
\end{theoreme}

Des cas particuliers de ce th{\'e}or{\`e}me ont {\'e}t{\'e} d{\'e}montr{\'e}s
dans \cite{GGR}, \cite{GO1}, \cite{GO2}, \cite{GM} ainsi que dans
\cite{Roberts}. Une d{\'e}monstration g{\'e}n{\'e}rale se trouve dans
\cite{truval}, mais selon Charles Walter, elle contient un point litigieux au
milieu de la page 116 (`` ...it is clear that the monomial $\gamma$ is greater
than any other $M-$product ...''). La d{\'e}monstration donn{\'e}e ici est
plus de nature  g{\'e}om{\'e}trique qu'alg{\'e}brique et est un raffinement
de la m{\'e}thode d'Horace introduite par A. Hirschowitz et C. Simpson dans
\cite{hirsim} pour prouver la ``Minimal Resolution Conjecture'' (MRC) de
\cite{mrc} pour tout nombre de points suffisamment grand et en position
suffisamment g{\'e}n{\'e}rale dans $\pp^n$.

\paragraphe{R{\'e}duction du probl{\`e}me}

Pour {\'e}tablir le th{\'e}or{\`e}me~\ref{th1} il suffit {\'e}videment de
prouver le suivant.

\begin{letheoreme}
  Pour tout entier non-n{\'e}gatif $a$ il existe des points $P_1,\dots,P_a$ de
  $\pp^n$, tels que tout entier $\ell$, le morphisme d'{\'e}valuation
  $$
    \sigma_{\ell,a} : H^0(\pp^n,\tp{n}(\ell)) \rightarrow
    \tp{n}(\ell)|_{P_1}\oplus\ldots\oplus \tp{n}(\ell)|_{P_a}
  $$
  est de rang maximal.
\end{letheoreme}
(on utilise l'isomorphisme $\Omega_{\pp^n}^{n-1}(\ell) \simeq
T_{\pp^n}(\ell-n-1)$.)\vspace{5mm}

On va d'abord r{\'e}duire le probl{\`e}me de rang maximal {\`a} d{\'e}calage
$\ell$ fix{\'e} {\`a}
celui d'une bijection, ceci gr{\^a}ce {\`a} des arguments standards (cf
\cite{harthirsch}, \cite{Ballico}, \cite{thIda}).
Pour plus de simplicit{\'e}, on posera

\begin{itemize}
\itemb $t_n(\ell) := h^0(\pp^n,\tp{n}(\ell))$,
\itemb $o_n(\ell) := h^0(\pp^n,\op{}{n}(\ell))$.
\end{itemize}

\begin{lelemme}
  Soient $q = q(\ell)$ et $r = r(\ell)$ respectivement le quotient et le
  reste de la {division eu-} clidienne
  de $t_n(\ell)$ par $n$. Alors pour que le th{\'e}or{\`e}me 1 soit vrai, il
  suffit qu'il existe un point $P_{q+1}$ tel que pour tout quotient
  $$
    \tp{n}(\ell)_{P_{q+1}} \rightarrow B \rightarrow 0
  $$
  de dimension $r$, il existe des points $P_1,\dots,P_{q+1}$ tels que
  le morphisme d'{\'e}valuation
  $$
   \tau_\ell : H^0(\pp^n,\tp{n}(\ell)) \rightarrow \tp{n}(\ell)|_{P_1}
   \oplus\ldots\oplus\tp{n}(\ell)|_{P_q}\oplus B
  $$
  soit bijectif.
\end{lelemme}
\demo
Supposons donc qu'il existe des points $P_1,\dots,P_{q+1}$ tels que
$\tau_\ell$ soit bijectif. Alors si $z \leq q$, on a une application
$\sigma_{\ell,z}$, compos{\'e}e de $\tau_{\ell}$ et de la surjection
$$
  \begin{array}{c}
    \tp{n}(\ell)|_{P_1}\oplus\ldots\oplus \tp{n}(\ell)|_{P_q}\oplus B \\
                  \downarrow \\
    \tp{n}(\ell)|_{P_1}\oplus\ldots\oplus \tp{n}(\ell)|_{P_z}
  \end{array}
$$
donc surjective.

Soit maintenant $z'$ un entier plus grand que $q+1$. Soient alors
$P_{q+2},\dots P_{z'}$ des points distincts, et diff{\'e}rents de
$P_1,\dots,P_{q+1}$. A l'aide du quotient $\tp{n}(\ell)_{P_{q+1}}\rightarrow
B \rightarrow 0$, on d{\'e}duit que le morphisme $\tau_\ell$
se factorise par $\sigma_{\ell,z'}$ et la projection
$$
  \begin{array}{c}
    \tp{n}(\ell)|_{P_1}\oplus\ldots\oplus \tp{n}(\ell)|_{P_{z'}}\\
                  \downarrow \\
    \tp{n}(\ell)|_{P_1}\oplus\ldots\oplus \tp{n}(\ell)|_{P_q} \oplus B.
  \end{array}
$$
On en d{\'e}duit donc que $\sigma_{\ell,z'}$ est injective. \eop

Pour obtenir le r{\'e}sultat escompt{\'e} pour $a_{n-2}$ et $b_{n-1}$ il
suffit donc
de montrer le
\begin{letheoreme}
\label{thtpn}
  Pour tout entier $\ell$, il existe des points $P_1,\dots P_{q+1}$ tels que
  le morphisme $\tau_\ell$ soit bijectif.
\end{letheoreme}

Remarquons que le th{\'e}or{\`e}me est trivialement vrai si $\ell \leq -2$. on
fera donc dans la suite l'hypoth{\`e}se que $\ell \geq -1$.\newpage

\paragraphe{M{\'e}thodes d'Horace vectorielles}

Dans ce paragraphe on introduit la m{\'e}thode d'Horace pour les probl{\`e}mes
d'{\'e}valuation de sections de fibr{\'e}s vectoriels en des points et des
``fractions de points''. Cette m{\'e}thode,  essentiellement bas{\'e}e sur les
transformations {\'e}l{\'e}mentaires de fibr{\'e}s vectoriels, fut introduite
en 1984 dans une lettre d'A. Hirschowitz {\`a} R. Hartshorne pour montrer que
si  $P_1,\dots,P_{28}$ sont des points en position g{\'e}n{\'e}rale dans
$\pp^3$, l'application naturelle
$$
  H^0(\pp^3,\Omega_{\pp^3}(5)) \rightarrow \Omega_{\pp^3}(5)|_{P_1}
  \oplus \dots \oplus  \Omega_{\pp^3}(5)|_{P_{28}}
$$
est bijective. Elle a aussi {\'e}t{\'e} utilis{\'e}e par M. Id{\`a} dans
\cite{thIda} cette fois-ci avec des droites et des points pour calculer la
r{\'e}solution minimale des id{\'e}aux d'arrangement de droites en position
g{\'e}n{\'e}rale dans $\pp^3$, et par O.F. Rahavandrainy dans~\cite{felix}
pour calculer des r{\'e}solutions de fibr{\'e}s instantons. La
pr{\'e}sentation suivante est celle de A. Hirschowitz et C. Simpson
(cf.~\cite{hirsim}).

Fixons quelques notations qu'on utilisera r{\'e}guli{\`e}rement par la suite.
Soit $X$ une vari{\'e}t{\'e} projective lisse, et $X'$ un diviseur
non-singulier de $X$. Soit $\Ff$ un faisceau localement libre sur $X$ et
$$
0 \rightarrow \Ff '' \rightarrow \Ff|_{X'} \rightarrow \Ff ' \rightarrow 0
$$
une suite exacte stricte de faisceaux localement libres sur $X'$. On notera
$\Ee$ le noyau du morphisme $\Ff \rightarrow \Ff '$; on notera que $\Ee$ est
localement libre sur $X$ et qu'on a la suite exacte
$$
0 \rightarrow \Ff '(-X') \rightarrow \Ee|_{X'} \rightarrow \Ff ''
\rightarrow 0.
$$
{~}\newline\newline

\begin{flushleft}
  \underline{La m{\'e}thode d'Horace "simple" :}
\end{flushleft}

Le lemme suivant est une cons{\'e}quence triviale du lemme du serpent.
\begin{lemme}
  \label{lemA}
  Supposons donn{\'e} une application lin{\'e}aire bijective $\lambda :
  H^0(X',\Ff ' )\rightarrow L.$ Supposons que $H^1(X,\Ee )=0$. Soit
  $\mu : H^0 (X, \Ff )\rightarrow M$ une application lin{\'e}aire. Alors,
  pour que le morphisme
  $$
    H^0(X,\Ff)\rightarrow M\oplus L
  $$
  soit de rang maximal, il faut et il suffit que le morphisme
  $$
    H^0(X,\Ee)\rightarrow M
  $$
  le soit.
\end{lemme}
\newpage

\begin{flushleft}
  \underline{La m{\'e}thode d'Horace "diff{\'e}rentielle":}
\end{flushleft}

Le lemme suivant est le lemme 1 de~\cite{hirsim}.
\begin{lemme}
\label{lemB}
  On se donne maintenant une  application lin{\'e}aire  surjective
  $\lambda :H^0(X',\Ff ' )\rightarrow L$ et supposons qu'il existe un point
  $Z'\in X'$ tel que l'application  $H^0 (X', \Ff ' )\rightarrow L
  \oplus \Ff ' _{Z'}$ soit injective.

    Supposons encore que $H^1(X,\Ee )=0$. Alors il existe un quotient
  $\Ee _{Z'}\rightarrow D$ avec noyau contenu dans $\Ff ' (-X')_{Z'}$,
  de dimension $dim (D) = r(\Ff ) -dim (ker (\lambda )),$
  ayant la propri{\'e}t{\'e} suivante.

     Soit $\mu : H^0 (X, \Ff )\rightarrow M$ une application lin{\'e}aire.
  Pour qu'il existe $Z\in X$ tel que l'application
  $$
    H^0 (X, \Ff )\rightarrow M \oplus L \oplus \Ff _Z
  $$
  soit de rang maximal, il suffit que
  $$
    H ^0 ( X, \Ee )\rightarrow M \oplus D
  $$
  le soit.
\end{lemme}~\newline

Consid{\'e}rons maintenant $\Gg$ un quotient localement libre de $\Ff$ et
$\Kk$ le noyau du morphisme $\Ff \rightarrow \Gg$. Posons $\Hh =
(\Ee+\Kk)/\Kk$ le sous-faisceau de $\Gg$ engendr{\'e} par $\Ee$ et $\Gg'$ le
quotient $\Gg/\Hh$. G{\'e}om{\'e}triquement on a
$$
  \pp(\Gg') = \pp(\Gg)\cap\pp(\Ff') \subset \pp(\Ff).
$$
On supposera $\Gg'$ localement libre sur $X'$ de sorte que $\Hh$ est
localement libre sur $X$.\newline
Le lemme suivant est alors une g{\'e}n{\'e}ralisation du pr{\'e}c{\'e}dent.

\begin{lemme}
  \label{lemC}
  Soit $\lambda : H^0(X',\Ff') \rightarrow L$ une application lin{\'e}aire
  surjective. On suppose qu'il existe $Y' \in X'$ tel que l'application
  $H^0(X',\Ff')\rightarrow L\oplus \Gg'_{Y'}$ soit injective. On suppose
  encore que $H^1(X,\Ee)=0$. Alors il existe un quotient
  $\Hh_{Y'}\rightarrow D$ avec noyau contenu dans $\Gg(-X')$ et de dimension
  $rg(\Gg)-dim(ker(\lambda))$ ayant la propri{\'e}t{\'e} suivante.

  Soit $\mu : H^0 (X, \Ff )\rightarrow M$ une application lin{\'e}aire. Pour
  qu'il existe $Y\in X$ tel que l'application
  $$
    H^0 (X, \Ff )\rightarrow M \oplus L \oplus \Gg_Y
  $$
  soit de rang maximal, il suffit que
  $$
   H ^0 ( X, \Ee )\rightarrow M \oplus D
  $$
  le soit.
\end{lemme}
\demo
Puisque l'application $H^0(X',\Ff')\rightarrow L\oplus \Gg'_{Y'}$ est
injective, l'application
$$
  H^0(X',\Ff')\rightarrow L\oplus \Ff'_{Y'}
$$
l'est aussi. Ces conditions ainsi que la premi{\`e}re hypoth{\`e}se du lemme
entra{\^\i}nent que le morphisme
$$
  ker(\lambda)\otimes_k\Oo_{X'} \rightarrow \Ff'\;\;\;(resp.\;
  ker(\lambda)\otimes_k\Oo_{X'} \rightarrow \Gg')
$$
est injectif et que son image $\Aa' \subset \Ff'$ (resp. $\Aa'' \subset
\Gg'$) est un sous-fibr{\'e} de $\Ff'$ (resp. $\Gg'$) au voisinage de $Y'$.
Notons que $rg(\Aa') = rg(\Aa'') = dim(ker(\lambda))$ et que $\Aa''$ est
l'image de $\Aa'$ par le morphisme surjectif $\Ff'\rightarrow \Gg'$.

Notons $\Bb$ le noyau de $\Ff\rightarrow \Ff'/\Aa'$. Alors le noyau de
$H^0(X,\Ff)\rightarrow L$ est {\'e}gal {\`a} $H^0(X,\Bb)$. Soit $\Cc$ l'image
de
$\Bb$ par le morphisme $\Ff \rightarrow \Gg$. Soit $D$ l'image du morphisme
 $\Hh_{Y'} \rightarrow \Cc_{Y'}$. On a une suite exacte
$$
  0\rightarrow D \rightarrow \Cc_{Y'} \rightarrow \Aa''_{Y'} \rightarrow 0.
$$
Puisque $\Aa'' \subset \Gg'$ est un sous-fibr{\'e} de rang $dim(ker(\lambda))$
 au voisinage de $Y'$, et $\Cc$ est localement libre au voisinage de $Y'$
avec $rg(\Cc) = rg(\Gg)$, on a
$$
  dim(D) = rg(\Gg) - dim(ker(\lambda)).
$$

Si $Y$ est un point de $X-X'$ les fibres $\Cc_Y$ et $\Gg_Y$ s'identifient.
Le fait que $H^1(X,\Ee)$ est nul entra{\^\i}ne que l'application
$H^0(X,\Ff) \rightarrow H^0(X',\Ff')$ est surjective, en particulier la
suite
$$
  0\rightarrow H^0(X,\Bb) \rightarrow H^0(X,\Ff) \rightarrow L \rightarrow 0
$$
est exacte, et puisque la suite
$$
  0 \rightarrow M\oplus\Cc_Y \rightarrow M\oplus L \oplus \Gg_Y
    \rightarrow L \rightarrow 0
$$
l'est aussi, on obtient alors que le morphisme
$$
  H^0(X,\Ff) \rightarrow M\oplus L \oplus \Gg_Y
$$
est de rang maximal  si et seulement si
$$
H^0(X,\Bb) \rightarrow M \oplus \Cc_Y
$$
l'est aussi.

On va maintenant sp{\'e}cialiser et prendre $Y = Y'$. Puisque $H^1(X,\Ee)=0$,
on a une suite exacte
$$
0\rightarrow H^0(X,\Ee) \rightarrow H^0(X,\Bb) \rightarrow H^0(X,\Aa')
 \rightarrow  0
$$
Du diagramme commutatif {\`a} lignes exactes
$$\diagram
0\rto & \Ee \rto\dto & \Bb \rto\dto& \Aa'\rto\dto & 0\\
0\rto & \Hh\rto      & \Cc \rto    & \Aa'' \rto   & 0
\enddiagram$$
on d{\'e}duit que le morphisme naturel compos{\'e}
$$
  H^0(X,\Ee) \subset H^0(X,\Bb) \rightarrow \Cc_{Y'} \rightarrow \Aa''_{Y'}
$$
est nul.
On en d{\'e}duit que le morphisme naturel compos{\'e}
$$
 H^0(X,\Bb) \rightarrow M\oplus\Cc_{Y'} \rightarrow \Aa''_{Y'}
$$
se factorise par $H^0(X,\Aa') \rightarrow \Aa''_{Y'}$.
Le morphisme $H^0(X,\Ee)\rightarrow M\oplus \Cc_{Y'}$ se factorise alors
par le noyau de  $M\oplus \Cc_{Y'} \rightarrow \Aa_{Y'}$, c'est {\`a} dire
$M \oplus D$. Le noyau du morphisme $H^0(X,\Bb) \rightarrow M\oplus \Cc_{Y'}$
est contenu dans le noyau de
$$
  H^0(X,\Ff)\rightarrow L \oplus \Aa''_{Y'} \subset L\oplus \Gg'_{Y'}.
$$
Par hypoth{\`e}se, il est donc contenu dans le noyau de
$H^0(X,\Ff)\rightarrow H^0(X,\Ff')$, c'est {\`a} dire $H^0(X,\Ee)$ et donc les
noyaux des deux morphismes suivants son {\'e}gaux :
$$
  H^0(X,\Ee) \rightarrow M\oplus D,\;\; H^0(X,\Bb)
  \rightarrow M\oplus \Cc_{Y'}
$$
Puisque le morphisme $H^0(X',\Aa') \rightarrow \Aa''_{Y'}$ est surjectif, on
conclut alors que pour que
$$
  H^0(X,\Bb) \rightarrow M\oplus \Cc_{Y'}
$$
soit de rang maximal, il suffit que
$$
  H^0(X,\Ee) \rightarrow M \oplus D
$$ le soit. Mais la condition de maximalit{\'e} du rang {\'e}tant ouverte, on
en d{\'e}duit alors qu'il existe $Y \in X-X'$ tel que
$$
  H^0(X,\Bb) \rightarrow M\oplus \Cc_Y
$$
soit de rang maximal, ce qui donne la conclusion cherch{\'e}e.\eop

\paragraphe{Les {\'e}nonc{\'e}s R, RB et un lemme d'Horace}

\subsubsection*{Enonc{\'e} ${\bf R}(\Ff;a)$}

Soit $a$ un entier non-n{\'e}gatif.  L'{\'e}nonc{\'e} ${\bf R}(\Ff ; a)$ veut
dire que pour tout entier non-n{\'e}gatif $z$, il existe des points
$T_1,\dots,T_a \in X$, tels que pour tous quotients
$$
  \Ff _{T_i}\rightarrow A_i \rightarrow 0,
$$
il existe des points $P_1,\dots,P_z \in X$, tels que
le morphisme d'{\'e}valuation
\begin{eqnarray*}
  H^0 (X,\Ff) &\rightarrow&  A_1 \oplus  \ldots \oplus A_a \oplus \\
  & & \Ff_{P_1}\oplus \dots \oplus\Ff_{P_z}
\end{eqnarray*}
 est de rang maximal.

\subsubsection*{Enonc{\'e} ${\bf RB}(\Ff,\Ff',z,y;\alpha,\beta)$}
Soit $X'$ un diviseur lisse connexe de $X$ et soient $z,y,\alpha$ et $\beta$
des entiers non-n{\'e}gatifs v{\'e}rifiant
\begin{enumerate}
\itemb $rg(\Ff)z + rg(\Ff')y + \alpha + \beta = h^0(X,\Ff)$
\itemb $rg(\Ff')y + \alpha + b \leq h^0(X,\Ff')$
\itemb $0\leq \alpha\leq rg(\Ff')$
\itemb si $\beta \not=0$, $rg(\Ff')\leq \beta < rg(\Ff)$,
\end{enumerate}
o{\`u} $b = b(\beta)$ vaut $r'$ si  $\beta$ est non nul et 0 sinon.

L'{\'e}nonc{\'e} ${\bf RB}(\Ff,\Ff',z,y;\alpha,\beta)$ signifie alors la chose
suivante. Il existe des points $T$ et $V$ de $X'$, tels que pour tout
quotient
$$
  \Ff'_T \rightarrow A \rightarrow 0
$$
de dimension $\alpha$ et si $\beta \not=0$,pour tout quotient d{\'e}pendant
rationnellement de $y$ points g{\'e}n{\'e}raux de $X'$,
$$
  \Ff_V\rightarrow B \rightarrow 0
$$
de dimension $\beta$ avec noyau contenu dans $\Ff''_V$, (si $\beta = 0$ on
pose $B= 0$), il existe $y$ points $Y_1,\dots,Y_y \in X'$ et $z$ points
$P_1,\dots,P_z \in X$ tel que le morphisme d'{\'e}valuation
\begin{eqnarray*}
  H^0(X,\Ff) &\rightarrow & A \oplus B \oplus \\
             &  & \Ff'_{Y_1}\oplus\dots\oplus \Ff'_{Y_y}\oplus \\
             &  & \Ff_{P_1}\oplus\dots\oplus \Ff_{P_z}
\end{eqnarray*}
soit bijectif.\newline\newline

Soient $z,y,\alpha$ et $\beta$ des entiers v{\'e}rifiant les conditions de
l'{\'e}nonc{\'e} ${\bf RB}(\Ff,\Ff',z,y;\alpha,\beta)$. On rappelle qu'on
avait d{\'e}fini l'entier $b=b(\beta)$ comme  valant $r'$ si
$\beta$ est non nul et 0 sinon. Posons alors
\begin{itemize}
\itemb $t := h^0(X',\Ff')- r'y-\alpha - b$,
\itemb $y'$ le quotient de la division euclidienne de $t$ par $r'$, $\delta$
       le reste,
\itemb $\zeta := 1$ si $\delta \not=0$ et $0$ sinon,
\itemb $\beta' := \zeta(r-\delta)$
\itemb $\alpha' := \beta - r'$ si $\beta \not= 0$ et 0 sinon.
\itemb $z' := z-y'-\zeta$ qu'on supposera non-n{\'e}gatif.
\end{itemize}
o{\`u} $r := rg(\Ff)=rg(\Ee)$ et $r' := rg(\Ff')$. On peut alors {\'e}noncer
le lemme suivant :
\begin{lemme}
  \label{lemRB}
  Soient $z,y,\alpha$ et $\beta$ des entiers non-n{\'e}gatifs v{\'e}rifiant
  les conditions de l'{\'e}nonc{\'e} ${\bf RB}(\Ff,\Ff',z,y;\alpha,\beta)$. On
  d{\'e}finit $a(\alpha) = 0$ si $\alpha$ est nul, 1 sinon. Supposons que
  $H^1(X,\Ee)=0$ et que (avec les notations pr{\'e}c{\'e}dant le lemme)
  les {\'e}nonc{\'e}s
$$
  {\bf R}(\Ff';a)\;et\;{\bf RB}(\Ee,\Ff'',z',y';\alpha',\beta')
$$
  soient vrais. Alors l'{\'e}nonc{\'e} ${\bf RB}(\Ff,\Ff',z,y;\alpha,\beta)$
  l'est aussi.
\end{lemme}
\demo
On ne va traiter que le cas o{\`u} $\beta \not= 0$, le cas $\beta=0$ {\'e}tant
laiss{\'e} au lecteur.
l'hypoth{\`e}se ${\bf R}(\Ff';a)$ entra{\^\i}ne qu'il existe $T$ et $V \in X'$,
tels que pour tout quotient
$$
  \Ff'_U\rightarrow A \rightarrow 0
$$ de dimension $\alpha$ et pour tout quotient
$$
  \Ff_V\rightarrow B \rightarrow 0,
$$
avec noyau contenu dans $\Ff''_V$
il existe des points $Y_1,\dots Y_y \in X'$ et $Z_1,\dots Z_{y'}$ tel que
le morphisme d'{\'e}valuation
\begin{eqnarray*}
  \lambda : H^0(X'\Ff') &\rightarrow& A\oplus \\
                        & & \Ff'_{Y_1}\oplus\dots\oplus\Ff'_{Y_y}\oplus\\
                        & & \Ff'_{Z_1}\oplus\dots\oplus\Ff'_{Z_{y'}}\\
                        & & \Ff'_V
\end{eqnarray*}
soit de rang maximal, donc surjectif ici.

$\bullet$ Si $\lambda$ est bijectif, (c'est en particulier le cas si
$rg(\Ff') = 1$), on entre dans le cadre du lemme~\ref{lemA} et on conclut
que pour tous $Z_{y'+1},\dots,Z_z \in X$,
le morphisme
\begin{eqnarray*}
  H^0(X,\Ff) &\rightarrow&   A\oplus B \oplus \\
             & & \Ff'_{Y_1}\oplus\dots\oplus\Ff'_{Y_y}\oplus \\
             & & \Ff_{Z_1}\oplus\dots\oplus\Ff_{Z_z}
\end{eqnarray*}
est bijectif pourvu que, en d{\'e}signant par $B''$ l'image de $\Ff''_V$ dans
$B$, le morphisme
\begin{eqnarray*}
  \epsilon : H^0(X,\Ee) &\rightarrow&  B'' \oplus\\
             & & \Ff''_{Z_1} \oplus\dots\oplus\Ff''_{Z_{y'}}\oplus \\
             & & \Ee_{Z_{y'+1}}\oplus\dots\oplus\Ee_{Z_z}
\end{eqnarray*}
le soit. L'hypoth{\`e}se ${\bf RB}(\Ee,\Ff'',z',y',\alpha',\beta')$ garantit
qu'il existe un point $V$  et des points $Z_1,\dots,Z_z$ tel que $\epsilon$
soit de rang maximal. L'existence de tels choix pour l'une ou l'autre des
hypoth{\`e}ses entra{\^\i}ne l'existence d'ouverts de choix qui sont
irr{\'e}ductibles,
puisqu'ils sont des ouverts de produit de $X$ ou $X'$. L'ouvert des choix de
la premi{\`e}re hypoth{\`e}se intersecte donc celui de la seconde.
On notera le fait important que $B''$ d{\'e}pend a priori de
$Y_1,\dots,Y_y$, lesquels {\em n'interviennent pas} dans $\epsilon$.\newline

$\bullet$ Si $\lambda$ n'est pas bijectif, il existe un point
$\bar{Z}$ de $X'$ tel que le morphisme
\begin{eqnarray*}
  H^0(X'\Ff') &\rightarrow& A\oplus \\
              & & \Ff'_{Y_1}\oplus\dots\oplus\Ff'_{Y_y}\oplus\\
              & & \Ff'_{Z_1}\oplus\dots\oplus\Ff'_{Z_{y'}}\oplus\\
              & & \Ff'_V\oplus \\
              & &\Ff'_{\bar{Z}}
\end{eqnarray*}
est injectif. On entre alors dans le cadre du lemme~\ref{lemB} et on conclut
qu'il existe un quotient $\Ee_{\bar{Z}}\rightarrow D$, de dimension $\beta'$,
avec noyau contenu dans $\Ff'_{\bar{Z}}$, avec la propri{\'e}t{\'e} suivante.
Pour tout ensemble de points $Z_{y'+2},\dots,Z_z$, et tout quotient $B$, il
existe $Z_{y'+1}$ dans $X$ tel que le morphisme
\begin{eqnarray*}
 H^0(X,\Ff) &\rightarrow& A\oplus B \oplus \\
            & & \Ff'_{Y_1}\oplus\dots\oplus \Ff'_{Y_y}\oplus\\
            & & \Ff_{Z_1}\oplus\dots\oplus \Ff_{Z_{y'}}\oplus\\
            & & \Ff_{Z_{y'+1}}\oplus \\
            & & \Ff_{Z_{y'+2}}\oplus\dots\oplus \Ff_{Z_z}
\end{eqnarray*}
est bijectif, pourvu que, en d{\'e}signant par $B''$ l'image de $\Ff''_V$ dans
$B$, le morphisme
\begin{eqnarray*}
  \epsilon : H^0(X,\Ee) &\rightarrow&  B'' \oplus D\oplus\\
             & & \Ff''_{Z_1} \oplus\dots\oplus\Ff''_{Z_{y'}}\oplus \\
             & & \Ee_{Z_{y'+2}}\oplus\dots\oplus\Ee_{Z_z}
\end{eqnarray*}
le soit. L'hypoth{\`e}se ${\bf RB}(\Ee,\Ff'',z',y',\alpha',\beta')$ et
l'argument donn{\'e} dans le cas pr{\'e}c{\'e}dent permettent alors de
conclure.\eop

\paragraphe{L'{\'e}nonc{\'e} MB et un autre lemme d'Horace}
\subsubsection*{Enonc{\'e} ${\bf MB}(\Ff,\Gg,z,y;a)$}
Soit maintenant
$$
\Ff\rightarrow\Gg\rightarrow 0
$$
une suite exacte stricte de faisceaux localement libres et $z,y$ et $a$
des entiers v{\'e}rifiant
\begin{enumerate}
  \itemb $rg(\Ff)z + rg(\Gg)y +a  = h^0(X,\Ff)$,
  \itemb $rg(\Gg)(z+y) + a \leq h^0(X',\Ff')$
  \itemb $a < rg(\Gg)$.
\end{enumerate}\vspace{3mm}

L'{\'e}nonc{\'e} ${\bf MB}(\Ff,\Gg,z,y;a)$ signifie alors la chose suivante :
Il existe un point $V$ tel que pour tout quotient
$$
  \Gg_T \rightarrow A \rightarrow 0
$$
de dimension $a$, il existe des points $Z_1,\dots,Z_z$ et des points
$Y_1,\dots,Y_y$ dans $X$ tels que
le morphisme d'{\'e}valuation
\begin{eqnarray*}
   H^0(X,\Ff) &\rightarrow&  A \oplus\\
              & & \Gg_{Y_1}\oplus\dots\oplus\Gg_{Y_y}\\
              & & \Ff_{Z_1}\oplus\dots\oplus\Ff_{Z_z}
\end{eqnarray*}
soit bijectif.

On va maintenant donner un lemme d'Horace permettant d'obtenir un
{\'e}nonc{\'e} {\bf RB} {\`a} partir d'un {\'e}nonc{\'e} de type {\bf R} et
d'un {\'e}nonc{\'e} de type {\bf MB}.

\begin{lemme}
  \label{lemMB}
  Soient $z,y$ et $a$ des entiers v{\'e}rifiant les conditions de
  l'{\'e}nonc{\'e}
  ${\bf RB}(\Ff,\Ff',z,y;a,0)$. Posons
  $$
    z' = \frac{h^0(X',\Ff|_{X'})-r'y-a}{r}
  $$
  qu'on suppose {\em entier} et  supposons que $H^1(X,\Ff(-X'))=0$ et que
  les {\'e}nonc{\'e}s
  $$
    {\bf R}(\Ff(-X');0)\;\;et\;\;{\bf MB}(\Ff|_{X'},\Ff',z',y,a)
  $$
  soient vrais. Alors l'{\'e}nonc{\'e} ${\bf RB}(\Ff,\Ff',z,y;a,0)$ l'est.
\end{lemme}
\demo
C'est une cons{\'e}quence triviale du lemme~\ref{lemA} et de la suite exacte
$$
  0\rightarrow\Ff(-X')\rightarrow\Ff\rightarrow\Ff_{|X'}\rightarrow 0.
$$ \eop

\paragraphe{Le th{\'e}or{\`e}me principal}
Le th{\'e}or{\`e}me suivant, contient dans sa partie (i) le th{\'e}or{\`e}me
\ref{thtpn} avec $z = q(\ell),y=a=0$ et $b = r(\ell)$.
\begin{theoreme}
  \label{thmprinc}
  Soit $n$ un entier $\geq 1$. Alors
  \begin{description}
  \item[(i)] on a ${\bf RB}(\tp{n}(\ell),\op{}{n-1}(\ell+1),z,y;0,b)$
             avec $y \leq o_{n-1}(\ell+1)$ et  $b \in \{0,\ldots,n-1\}$,
  \item[(ii)] on a ${\bf RB}(\op{n}{n}(\ell+1),\tp{n-1}(\ell),z,y;a,0)$
             avec $a \in \{0,\ldots, n-2\}$,
  \item[(iii)] on a ${\bf R}(\tp{n\!-\!1}(\ell);1)$,
  \item[(iv)] on a ${\bf MB}(\op{(n\!+\!1)}{n}(\ell+1),\tp{n}(\ell);z,y;a)$
             avec $z\geq o_n(\ell)$ et $a  \leq n-1$.
  \end{description}
\end{theoreme}

\begin{rem}
  Dans l'{\'e}nonc{\'e} (ii) on a toujours $z \geq o_n(\ell)$. En effet,des
  deux premi{\`e}res conditions des {\'e}nonc{\'e}s {\bf RB}  on tire que
  $$
    nz \geq no_n(\ell+1)-t_{n-1}(\ell) = t_{n}(\ell-1)
  $$
  et ce dernier est plus grand que $no_n(\ell)$.
\end{rem}

Le reste de cet article sera consacr{\'e} {\`a} la preuve de ce
th{\'e}or{\`e}me. Celle-ci se fera par r{\'e}currence sur l dimension. Pour
amorcer cette r{\'e}currence, notons d'abord que pour $n=1$, les assertions
(i) {\`a} (iii) sont trivialement vraies. On va maintenant donner une preuve
de l'assertion (iv) dans ce cas. Cette assertion  s'{\'e}crit ainsi.
\newline
On a
\begin{itemize}
  \itemb ${\bf MB}(\op{2}{1}(\ell+1),\op{}{1}(\ell+2),o_1(\ell+1),0;0)$,
  \itemb ${\bf MB}(\op{2}{1}(\ell+1),\op{}{1}(\ell+2),o_1(\ell),1;1)$.
\end{itemize}
Le premier de ces {\'e}nonc{\'e}s est trivialement vrai. Montrons alors le
second. Consid{\'e}rons donc des
points {\it distincts} $Z_1,\dots,Z_{o_1(\ell)},Y$ et $U$ de $\pp^1$. Il faut
alors montrer que le morphisme
\begin{eqnarray*}
  \mu : H^0(\pp^1,\op{2}{1}(\ell+1))&\rightarrow&
        \op{2}{1}(\ell+1)_{Z_1}\oplus\ldots\oplus
        \op{2}{1}(\ell+1)_{Z_{o_1(\ell)}}\oplus \\
        & & \op{}{1}(\ell+2)_Y\oplus \op{}{1}(\ell+2)_U
\end{eqnarray*}
est bijectif. On va utiliser pour ce faire la suite exacte suivante
$$
  0\rightarrow\op{}{1}(\ell)\rightarrow \op{2}{1}(\ell+1)\rightarrow
  \op{}{1}(\ell+2)\rightarrow 0,
$$
ainsi que celle qui s'en d{\'e}duit en cohomologie. Cette suite permet aussi,
si $Z\in \pp^1$ de d{\'e}composer $\op{2}{1}(\ell+1)_Z$ en
$\op{}{1}(\ell)_Z\oplus \op{}{1}(\ell+2)_Z$. Posons alors
\begin{eqnarray*}
  M &=& \op{}{1}(\ell)_{Z_1}\oplus\ldots\oplus
        \op{}{1}(\ell)_{Z_{o_1(\ell)}}\,,\\
  L &=& \op{}{1}(\ell+2)_{Z_1}\oplus\ldots\oplus
        \op{}{1}(\ell+2)_{Z_{o_1(\ell)}}\oplus \\
    & & \op{}{1}(\ell+2)_Y\oplus\op{}{1}(\ell+2)_U\,.
\end{eqnarray*}
Alors les applications $H^0(\pp^1,\op{}{1}(\ell))\rightarrow M$ et
$H^0(\pp^1,\op{}{1}(\ell+1))\rightarrow L$ sont bijectives. On en d{\'e}duit
alors que $\mu $ l'est. Le th{\'e}or{\`e}me est donc montr{\'e} pour $n=1$. On
supposera donc dans la suite $n\geq 2$.
\vspace{6mm}

On proc{\`e}de maintenant {\`a} la preuve des parties (i) et (ii)
du th{\'e}or{\`e}me. Elle se fait par r{\'e}ductions successives, pour aboutir
{\`a} un {\'e}nonc{\'e} du type
${\bf MB}(\op{n}{n-1}(m+1),\tp{n-1}(m),z,y;a:\alpha)$ avec
$z\geq o_{n-1}(m)$ et $\alpha \leq n-2$, qui est vrai par hypoth{\`e}se de
r{\'e}currence. On utilisera pour ce faire les trois lemmes de r{\'e}duction
suivants.
\newline

\begin{lemme}
  \label{lemred1}
  Soient $z,y$ et $b$ des entiers non n{\'e}gatifs v{\'e}rifiant les
  hypoth{\`e}ses de
  l'{\'e}nonc{\'e}
  $$
    {\bf  RB}(\tp{n}(\ell),\op{}{n}(\ell),z,y;0,b).
  $$
  Posons $\beta = 0$ si $b=0$, 1 sinon. Alors, pour que l'{\'e}nonc{\'e}
  ci-dessus
  soit vrai, il suffit que
  $$
    {\bf RB}(\op{n}{n}(\ell+1),\tp{n-1}(\ell),
    z-o_n(\ell+1)+y+\beta,o_n(\ell+1)-y-\beta;b-\beta,0)
  $$
  le soit.
\end{lemme}
\demo
C'est une instanciation du lemme~\ref{lemRB} (l'{\'e}nonc{\'e}
${\bf R}(\op{}{n\!-\!1}(\ell+1);1)$ est trivialement vrai). {\it Il faut
noter que le quotient intervenant dans l'{\'e}nou{\'e} r{\'e}duit est
ind{\'e}pendant des points intervenant dans cet {\'e}nonc{\'e}.} \newline
{~}\newline
Soient maintenant $z,y$ et $a$ des entiers non n{\'e}gatifs v{\'e}rifiant les
hypoth{\`e}ses de l'{\'e}nonc{\'e}
${\bf RB}(\op{n}{n}(\ell+1),\tp{n-1}(\ell),z,y;a,0)$.
Posons
\begin{itemize}
  \itemb $t = t_n(\ell)-3y-a$,
  \itemb $u$ le reste de la division euclidienne de $t$ par n-1,
  \itemb $b'$ = 0 si $u=0$, $b' = n-u$ sinon,
  \itemb $y'= (t-u)/(n-1)$ si $b' \not= 0$, $y= t/(n-1)$ sinon,
  \itemb $z'=z- y'-1$ si $b' \not= 0$, $z=z-y'$ sinon.
\end{itemize}

\begin{lemme}
  \label{lemred2}
  Supposons $t \leq n\,o_n(\ell)$.Alors pour que l'{\'e}nonc{\'e} ci-dessus
  soit vrai, il suffit que
  $$
    {\bf RB}(\tp{n}(\ell-1),\op{}{n}(\ell),z',y';0,b')
  $$
  le soit.
\end{lemme}
\demo
C'est une instanciation du lemme~\ref{lemRB} (on utilise l'hypoth{\`e}se de
r{\'e}currence sur ${\bf R}(\tp{n\!-\!1}(\ell);1))$.\newline

Supposons maintenant que $t > n\,o_n(\ell)$. On r{\'e}{\'e}crit cette
derni{\`e}re in{\'e}galit{\'e} en tenant compte de la relation $nz+(n-1)y+a =
no_n(\ell+1)$. On trouve que $z > o_n(\ell) + o_{n-1}(\ell)$. On a alors
le lemme suivant.
\begin{lemme}
  \label{alakon}
  Pour que l'{\'e}nonc{\'e} ${\bf RB}(\op{n}{n}(\ell+1),\tp{n}(\ell),z,y;a,0)$
  soit vrai, il suffit que
  l'{\'e}nonc{\'e}
  $$
    {\bf MB}(\op{n}{n\!-\!1}(\ell+1),\tp{n\!-\!1}(\ell),z-o_n(\ell),y;a)
  $$
  le soit.
\end{lemme}
\demo
C'est une instanciation du lemme~\ref{lemMB}.\newline

Pour d{\'e}montrer (ii) on utilise soit le lemme~\ref{alakon} si
l'{\'e}nonc{\'e} entre dans ses hypoth{\`e}ses, soit successivement les
lemmes~\ref{lemred2} et \ref{lemred1} autant qu'il est n{\'e}cessaire pour
obtenir soit un {\'e}nonc{\'e} trivialement vrai, soit rentrer dans les
hypoth{\`e}ses du lemme~\ref{alakon}. On utilise alors l'hypoth{\`e}se de
r{\'e}currence sur les {\'e}nonc{\'e}s {\bf MB} pour conclure.\eop
\vspace{6mm}

Pour d{\'e}montrer (i) on utilise le lemme~\ref{lemred1} qui permet alors de
ce r{\'e}duire {\`a}~(ii) et on conclut comme pr{\'e}c{\'e}demment.\eop
{~}\vspace{5mm}

On va maintenant proc{\'e}der {\`a} la d{\'e}monstration de la partie (iii)
du th{\'e}or{\`e}me. C'est une cons{\'e}quence  de
l' {\'e}nonc{\'e}
$$
  {\bf RB}(\tp{n}(\ell),\op{}{n-1}(\ell+1),q(\ell),0;0,r(\ell)).
$$

L'{\'e}nonc{\'e} ci-dessus  dit que pour tout $P_{q+1}\in \pp^{n-1}$,et
pour tout quotient
$$
  \tp{n}(\ell)_{Z_{q(\ell)+1}} \rightarrow B \rightarrow 0
$$
de dimension $r(\ell)$, il existe $P_1,\dots,P_q \in \pp^n$ tels que
le morphisme d'{\'e}valuation
$$
  H^0(\pp^n,\tp{n}(\ell)) \rightarrow \tp{n}(\ell)_{Z_1}
  \oplus\ldots\oplus \tp{n}(\ell)_{Z_{q(\ell)}}\oplus B
$$
soit bijectif.
Il suffira de montrer l'assertion pour $q(\ell)$ et $q(\ell)+1$ points
lorsque $a \leq r(\ell)$ et pour $q(\ell)-1$ et $q(\ell)$ points lorsque $a
> r(\ell)$.
Soit $\tp{n}(\ell)_{Z_{q+1}}\rightarrow A\rightarrow 0$ un quotient.
de dimension $a$. Supposons $a \leq r(\ell)$. Il existe alors un
quotient $\tp{n}(\ell)_{Z_{q+1}}\rightarrow B\rightarrow 0$ tel que $A$ soit un
quotient de $B$. On conclut alors d'apr{\`e}s ce qu'on a vu plus haut qu'il
existe des points $Z_1,\dots Z_{q(\ell)}$ tels que, en notant $L$ l'espace
vectoriel
$$
\tp{n}(\ell)_{Z_1} \oplus\ldots\oplus \tp{n}(\ell)_{Z_q}
$$
l'application
$$
  H^0(\pp^n,\tp{n}(\ell)) \rightarrow L\oplus A
$$
soit surjective.
Par semi-continuit{\'e} on d{\'e}duit qu'il existe un point $T \not=
Z_{q+1}$ et un quotient $\tp{n}(\ell)_T\rightarrow Q \rightarrow 0$ de
dimension $r(\ell)$ tel que l'application
$$
  H^0(\pp^n,\tp{n}(\ell)) \rightarrow L\oplus Q
$$
soit bijective.
De la surjection $L\oplus A \oplus \tp{n}(\ell)_Q\rightarrow L\oplus Q$ on
d{\'e}duit l'injectivit{\'e} de
$$
  H^0(\pp^n,\tp{n}(\ell)) \rightarrow L\oplus\tp{n}(\ell)_Q\oplus A.
$$
Le cas $a >r(\ell)$ se traite avec des arguments similaires et est
laiss{\'e} au lecteur. \eop

{~}\vspace{1cm}

Il reste maintenant, pour conclure la d{\'e}monstration du th{\'e}or{\`e}me
{\`a} prouver la partie (iv). Pour cela, on distinguera, dans l'{\'e}nonc{\'e}
$$
  {\bf MB}(\op{(n\!+\!1)}{n}(\ell+1),\tp{n}(\ell);z,y;a)
$$
quatre cas : $z = o_n(\ell+1)$, qui est trivialement vrai,
$o_n(\ell-1)+o_{n-1}(\ell+1) \leq z < o_n(\ell+1)$, qu'on r{\'e}duira
trivialement {\`a} un autre {\'e}nonc{\'e} du m{\^e}me type mais de degr{\'e}
plus petit, le cas $o_n(\ell)<z < o_n(\ell-1)+o_{n-1}(\ell+1)$ et enfin le
cas o{\`u} $z=o_n(\ell)$. On va donc prouver les trois cas restants.
\newline
\newline
\underline{Second cas :} $o_n(\ell-1)+o_{n-1}(\ell+1) \leq z < o_n(\ell+1)$.
Choisissons un hyperplan $\pp^{n-1} \subset \pp^n$ et soient $Z_1,\dots,
Z_{o_{n-1}(\ell+1)}\in \pp^{n-1}$ en position g{\'e}n{\'e}rale.
L'application d'{\'e}valuation des sections
$$
\begin{array}{c}
  H^0(\pp^{n-1},\op{(n+1)}{n}(\ell+1)_{|\pp^{n-1}})\\
                   \downarrow \\
  \op{(n\!+\!1)}{n}(\ell+1)_{|\pp^{n-1}\,Z_1}\oplus\dots\oplus
  \op{(n\!+\!1)}{n}(\ell+1)_{|\pp^{n-1}\,Z_{o_{n-1}(\ell+1)}}
\end{array}
$$
est alors bijective. Gr{\^a}ce au lemme~\ref{lemA}, avec la suite exacte
$$
  0\rightarrow \op{(n\!+\!1)}{n}(\ell)\rightarrow\op{(n\!+\!1)}{n}(\ell+1)
  \rightarrow \op{(n\!+\!1)}{n}(\ell+1)_{|\pp^{n-1}}\rightarrow 0
$$
on se r{\'e}duit alors
{\`a} l'{\'e}nonc{\'e}
$$
  {\bf MB}(\op{(n\!+\!1)}{n}(\ell),\tp{n}(\ell-1);z-o_{n-1}(\ell+1),y;a).
$$
\newline
\newline
\underline{Troisi{\`e}me cas :} $z=o_n(\ell)$. On utilisera pour cela la suite
d'Euler pour le fibr{\'e} tangent et celle qui s'en d{\'e}duit en cohomologie :
$$
\begin{array}{c}
  0\rightarrow \op{}{n}(\ell) \rightarrow \op{(n+1)}{n}(\ell+1)
  \rightarrow \tp{n}(\ell)\rightarrow 0 ,\\
  0\rightarrow H^0(\pp^n,\op{}{n}(\ell))\rightarrow
  H^0(\pp^n,\op{(n+1)}{n}(\ell+1)) \rightarrow
  H^0(\pp^n,\tp{n}(\ell)) \rightarrow 0.
\end{array}
$$

Remarquons que Si $P \in \pp^n$, alors, gr{\^a}ce {\`a} la suite d'Euler,
$\op{(\!n\!+\!1\!)}{n}(\!\ell\!+\!1)_P$ se d{\'e}compose en
$\op{}{n}(\ell)_P\oplus\tp{n}(\ell)_P$.
Les entiers $o_n(\ell)+y$ et $a$ v{\'e}rifient les hypoth{\`e}ses de
l'{\'e}nonc{\'e}
${\bf RB}(\tp{n}(\ell),\op{}{n-1}(\ell+1),,0;0,a)$.
Soit donc $V$ un point de $\pp^{n-1}$ tel que pour tout quotient
$$
  \tp{n}(\ell)_V \rightarrow A \rightarrow 0
$$
il existe des points $Z_1,\dots,Z_{o_n(\ell)},Y_1,\dots,Y_y \in \pp^n$
tels que en notant $L$ l'espace vectoriel
$$
  A\oplus\tp{n}(\ell)_{Z_1}\oplus\dots\oplus\tp{n}(\ell)_{Z_{o_n(\ell)}}\oplus
  \tp{n}(\ell)_{Y_1}\oplus\dots\oplus\tp{n}(\ell)_{Y_y}
$$
l'application $H^0(\pp^n,\tp{n}(\ell))\rightarrow L$ soit bijective.
Si $P_1,\dots,P_{o_n(\ell)} \in \pp^n$ sont en position suffisamment
g{\'e}n{\'e}rale, en notant $M$ l'espace vectoriel
$$
  \op{}{n}(\ell)_{P_1}\oplus\dots\oplus\op{}{n}(\ell)_{P_{o_n(\ell)}}
$$
l'application $H^0(\pp^n,\op{}{n}(\ell))\rightarrow M$ est bijective. Or
les ouverts de choix des $Z_i$ et des $P_j$ sont des ouverts de
$(\pp^n)^{o_n(\ell)}$ et donc s'intersectent. On peut donc supposer $P_i =
  Z_i$ et on conclut alors par le lemme du serpent que l'application
$$
 H^0(\pp^n,\op{(n+1)}{n}(\ell+1))\rightarrow M\oplus L
$$
est bijective.
\newline
\newline
\underline{Quatri{\`e}me cas :} $o_n(\ell)<z<o_n(\ell-1)+o_{n-1}(\ell+1)$.
C'est le plus technique d'entre eux et c'est pour les besoins cette
d{\'e}monstration qu'on a introduit le lemme~\ref{lemC}. On va d'abord
introduire les faisceaux et les suites exactes utilis{\'e}es dans la
d{\'e}monstration.

On a la suite exacte
$$
  0\rightarrow
  \op{n}{n}(\ell+1)\oplus\op{}{n}(\ell)\rightarrow
  \op{(n+1)}{n}(\ell+1)\rightarrow\op{}{n-1}(\ell+1)\rightarrow 0
$$
et par transformation {\'e}l{\'e}mentaire
$$
  0\rightarrow \op{(n+1)}{n}(\ell)\rightarrow \op{n}{n}(\ell+1)
  \oplus\op{}{n}(\ell)\rightarrow \op{n}{n-1}(\ell+1)\rightarrow 0.
$$

On a aussi le diagramme commutatif {\`a} lignes et colonnes exactes
$$\diagram
  0\rto& \op{n}{n}(\ell\!+\!1)\oplus\op{}{n}(\ell)\dto\rto&
  \op{(n\!+\!1)}{n}(\ell\!+\!1)\dto\rto&\op{}{n-1}(\ell\!+\!1)\ddouble\rto&
  0\\
  0\rto& \op{n}{n}(\ell\!+\!1)\rto\dto& \tp{n}(\ell)\dto \rto &
  \op{}{n-1}(\ell\!+\!1)\rto& 0\\  & 0 &  0 &
\enddiagram$$
et en appliquant la transformation {\'e}l{\'e}mentaire pr{\'e}c{\'e}dente, le
diagramme suivant commutatif {\`a} lignes et colonnes exactes.
$$\diagram
  0\rto& \op{(n+1)}{n}(\ell)\dto\rto&\op{n}{n}(\ell+1)\oplus
  \op{}{n}(\ell)\rto\dto&\op{n}{n-1}(\ell+1)\rto\dto&0\\
  0\rto& \tp{n}(\ell-1)\rto\dto & \op{n}{n}(\ell+1) \rto\dto &
  \tp{n-1}(\ell)\rto\dto&0\\
  & 0 & 0 & 0 &
\enddiagram$$

On posera $\alpha = 0$ si $a = 0$ et $\alpha = 1$ sinon. On notera encore
$\pp^{n-1}$ un hyperplan de $\pp^n$. Soit $d$ l'entier $z - o_n(\ell)$ et
$y'$ l'entier $o_{n-1}(\ell+1)-d-\alpha$. Les hypoth{\`e}ses entra{\^\i}nent
que $d$ et $y'$ sont tous deux non-n{\'e}gatifs.
Supposons $\alpha\not=0$. Alors il existe $V\in \pp^{n-1}$ tel que pour tout
quotient
$$
  \tp{n}(\ell)_V\rightarrow A \rightarrow 0
$$
de dimension $a$ avec noyau contenu dans $\tp{n-1}(\ell)$, il existe des
points $Z_1,\dots,Z_d, Y_1,\dots,Y_{y'} \in \pp^{n-1}$ tel que l'application
\begin{eqnarray*}
  H^0(\pp^{n-1},\op{}{n-1}(\ell+1)) &\rightarrow& \op{}{n-1}(\ell+1)_V
  \oplus \\
  & &
  \op{}{n-1}(\ell+1)_{Z_1}\oplus\dots\oplus\op{}{n-1}(\ell+1)_{Z_d}\oplus \\
  & &
  \op{}{n-1}(\ell+1)_{Y_1}\oplus\dots\oplus\op{}{n-1}(\ell+1)_{Y_{y'}}
\end{eqnarray*}
soit bijective et si $\alpha = 0$, l'application
\begin{eqnarray*}
  H^0(\pp^{n-1},\op{}{n-1}(\ell+1)) &\rightarrow&
  \op{}{n-1}(\ell+1)_{Z_1}\oplus\dots\oplus\op{}{n-1}(\ell+1)_{Z_d}\oplus \\
  & &
  \op{}{n-1}(\ell+1)_{Y_1}\oplus\dots\oplus\op{}{n-1}(\ell+1)_{Y_{y'}}
\end{eqnarray*}
l'est dans ce cas.
On rentre alors dans le cadre du lemme~\ref{lemA}et on conclut que,
pour que l'{\'e}nonc{\'e}
${\bf MB}(\op{(n\!+\!1)}{n}(\ell\!+\!1),\tp{n}(\ell),z,y;a)$ soit vrai, il
suffit qu'il existe des points $Z_{d+1},\dots,Z_z$, $ Y_{y'+1},\dots,Y_y \in
\pp^n$,  tel que, en d{\'e}signant par $A'$ l'image de $\tp{n-1}(\ell)_V$ dans
$A$ si $A\not=0$, $A'=0$ sinon, que l'application
\begin{eqnarray*}
\label{monmu}
  \mu :H^0(\pp^n,\op{(\!n\!+\!1\!)}{n}(\!\ell\!+\!1\!))
  &\rightarrow & A' \oplus \\
  & & \tp{n-1}(\ell)_{Y_1}\oplus\dots\oplus\tp{n-1}(\ell)_{Y_{y'}}\oplus \\
  & &
  \op{n}{n-1}(\ell+1)_{Z_1}\oplus\dots\oplus\op{n}{n-1}(\ell+1)_{Z_d}\oplus\\
  & &
  \op{n}{n}(\ell+1)_{Y_{y'+1}}\oplus\dots\oplus\op{n}{n}(\ell+1)_{Y_y}\oplus
  \\
  &&\op{n}{n}(\ell+1)\oplus\op{}{n}(\ell)_{Z_{d+1}}\oplus\dots\oplus
    \op{n}{n}(\ell+1)\oplus\op{}{n}(\ell)_{Z_{o_n(\ell)}}
\end{eqnarray*}
soit bijective.

Notons que $a' := dim(A') = 0$ si $a \leq 1$, $a-1$ sinon.
Soit alors $e$ l'entier $no_{n-1}(\ell+1) - dn - (n-1)y'-a'$, $f$ et $g$
respectivement le quotient et le reste de la division euclidienne de $e$ par
$n-1$. On va alors diviser la preuve de la bijectivit{\'e} de $\mu $ en deux
cas, le plus simple {\'e}tant celui o{\`u} $g=0$, le second, $g\not= 0$
utilisera le lemme 3. \newline

Dans le premier cas, le lecteur v{\'e}rifiera que les entiers $d,y'+f$ et $a'$
v{\'e}rifient les hypoth{\`e}ses de l'{\'e}nonc{\'e}
${\bf MB}(\op{n}{n-1}(\ell+1),\tp{n-1}(\ell),d,y'+f,a')$, qui est suppos{\'e}
vrai par hypoth{\`e}se de r{\'e}currence.

Il existe donc un point $V \in \pp^{n-1}$ tel que pour tout quotient
$$
 \tp{n-1}(\ell)_V \rightarrow A' \rightarrow 0
$$
il existe des points $Z_1,\dots,Z_d, Y_1,\dots ,Y_{y'+f}$ tel que
l'application
\begin{eqnarray*}
  H^0(\pp^{n-1},\op{n}{n-1}(\ell+1)) & \rightarrow & A'\oplus \\
  & & \tp{n-1}(\ell)_{Y_1}\oplus \dots\oplus
  \tp{n-1}(\ell)_{Y_{y'+f}}\oplus\\
  & & \op{n}{n-1}(\ell+1)_{Z_1}\oplus\dots\oplus\op{n}{n-1}(\ell+1)_{Z_d}
\end{eqnarray*}
soit bijective.
On rentre alors dans le cadre du lemme~\ref{lemA} et on en conclut qu'il
existe des points $Y_{y'+f+1},\dots,Y_y,Z_{d+1},\dots,Z_z \in
\pp^n$ tel que  $\mu $ soit bijective pourvu que l'application
\begin{eqnarray*}
  \epsilon: H^0(\pp^n,\op{(n+1)}{n}(\ell))&\rightarrow&
\op{}{n-1}(\ell)_{Y_{y'+1}}\oplus\dots\oplus\op{}{n-1}(\ell)_{Y_{y'+f}}\oplus
\\
&&\tp{n}(\ell-1)_{Y_{y'+f+1}}\oplus\dots\oplus\tp{n}(\ell-1)_{Y_y}\oplus\\
&&\op{(n+1)}{n}(\ell)_{Z_{d+1}}\oplus\dots\oplus\op{(n+1)}{n}(\ell)_{Z_z}
\end{eqnarray*}
le soit. Remarquons alors que $z-d=o_n(\ell-1)$. En utilisant alors un
argument similaire {\`a} celui de la preuve du troisi{\`e}me cas de (iv),
pour que $\epsilon$ soit
bijective, il suffit que l'application
\begin{eqnarray*}
  H^0(\pp^n,\tp{n}(\ell-1)))&\rightarrow&
\op{}{n-1}(\ell)_{Y_{y'+1}}\oplus\dots\oplus\op{}{n-1}(\ell)_{Y_{y'+f}}\oplus
\\
&&\tp{n}(\ell-1)_{Y_{y'+f+1}}\oplus\dots\oplus\tp{n}(\ell-1)_{Y_y}\oplus\\
&&\tp{n}(\ell-1)_{Z_{d+1}}\oplus\dots\oplus\tp{n}(\ell-1)_{Z_z}
\end{eqnarray*}
le soit. Or l'existence de choix de points pour que ce dernier
{\'e}nonc{\'e}  soit v{\'e}rifi{\'e} est garanti par l'{\'e}nonc{\'e}
${\bf RB}(\tp{n}(\ell-1), \op{}{n-1}(\ell),o_n(\ell)-1+y-y'-f,f;0,0)$.
Pour conclure, on utilise alors
l'irreductibilit{\'e} des espaces de param{\`e}tres pour garantir la
s{\'e}rie de choix que l'on vient de faire.\newline

Reste maintenant {\`a} traiter le cas o{\`u} $g$ est non nul.
On va montrer que si $V\in \pp^n$, alors pour tout quotient
$$
\tp{n}(\ell)\rightarrow A'\rightarrow 0
$$
il existe des points $Z_1,\dots,Z_d$, $Y_1,\dots,Y_{y'+f}$ et un point
$\bar{Y}$ in $\pp^{n-1}$ tel que en notant $L$ l'espace vectoriel
$$
 \tp{n-1}(\ell)_{Y_1}\oplus \dots\oplus
  \tp{n-1}(\ell)_{Y_{y'+f}}\oplus
 \op{n}{n-1}(\ell+1)_{Z_1}\oplus\dots\oplus\op{n}{n-1}(\ell+1)_{Z_d}
$$
l'application
$
\lambda : H^0(\pp^{n-1},\op{n}{n-1}(\ell+1))\rightarrow L\oplus A'
$
est surjective et que l'application d{\'e}duite
$
\lambda' : H^0(\pp^{n-1},\op{n}{n-1}(\ell+1))\rightarrow L\oplus A'\oplus
\tp{n-1}(\ell)_{\bar{Y}}
$
est injective. Pour ce faire on va alors distinguer
suivant la valeur de $g+a'$.
\begin{itemize}
\itemb $g+a' \leq n-1$. Alors les entiers $d,y'+f$ et $g+a'$ v{\'e}rifient les
conditions de l'{\'e}nonc{\'e}
$$
  {\bf MB}(\op{n}{n-1}(\ell+1),\tp{n-1}(\ell),d,y'+f,g+a').
$$
\itemb $g+a'>n-1$. Alors les entiers $d,y'+f+1$ et $g+a'-(n-1)$ v{\'e}rifient
les
conditions de l'{\'e}nonc{\'e}
$$
  {\bf MB}(\op{n}{n-1}(\ell+1),\tp{n-1}(\ell),d,y'+f+1,g+a'-(n-1)).
$$
\end{itemize}
Ces deux {\'e}nonc{\'e}s sont suppos{\'e}s vrais par hypoth{\`e}se de
r{\'e}currence. Consid{\'e}rons le cas $g+a'\leq n-1$. D'apr{\`e}s
l'hypoth{\`e}se correspondante, il existe $T\in \pp^{n-1}$ tel que pour tout
quotient
$$
  \tp{n-1}(\ell)_T\rightarrow G \rightarrow 0
$$
de dimension $g+a'$, il existe des points $Z_1,\dots,Z_d$,
$Y_1,\dots,Y_{y'+f} \in \pp^{n-1}$ tel que, en notant $N$ l'espace vectoriel
$$
  \tp{n-1}(\ell)_{Y_1}\oplus\dots\oplus\tp{n-1}(\ell)_{Y_{y'+f}}\oplus
  \op{n}{n-1}(\ell+1)_{Z_1}\oplus\dots\oplus\op{n}{n-1}(\ell+1)_{Z_d}
$$
l'application
$$
  \xi : H^0(\pp^{n-1},\op{n}{n-1}(\ell)) \rightarrow  G\oplus N
$$
soit bijective.
On pose alors $V=T$. Tout quotient $\tp{n-1}(\ell)\rightarrow A'$ se
factorise par un quotient $G$ de dimension $g+a'$. On en d{\'e}duit alors
la surjectivite de $\lambda$, avec $L=N$.
Montrons l'injectivit{\'e} de $\lambda'$.
Par semi-continuit{\'e}, il existe un ouvert $U \subset \pp^{n-1}$, un point
$\bar{Y} \not= T$ et un quotient $\tp{n-1}(\ell)_{\bar{Y}}\rightarrow G'$ de
dimension $g'+ a$ tel que l'application
$ H^0(\pp^{n-1},\op{n}{n-1}(\ell)) \rightarrow  G'\oplus N$ soit aussi
bijective. De la surjection $N\oplus A'\oplus \tp{n-1}(\ell)_{\bar{Y}}
\rightarrow N\oplus G'$ on tire l'injectivit{\'e} de $\lambda'$.\newline
Le cas $g+a'>n-1$ se traite de fa{\c{c}}on similaire, il suffit seulement de
remarquer qu'on a {\'e}videment $g+a'-(n-1) \leq a'$.\newline

On rentre alors dans le cadre du lemme~\ref{lemC} et on conclut qu'il existe
un quotient
$$
\tp{n-1}(\ell)_{\bar{Y}}\rightarrow D \rightarrow 0
$$
de dimension $\delta=n-dim\,Ker\lambda$ avec noyau contenu dans
$\tp{n-1}(\ell-1)_{\bar{Y}}$ jouissant de la propri{\'e}t{\'e} suivante.
Pour tout choix de points $Z_{d+1},\dots,Z_z$, $Y_{y'+f+2},\dots,Y_y \in
\pp^n$, il existe un point $Y_{y'+f +1} \in \pp^n$ tel que l'application
$\mu $ correspondante (page \pageref{monmu}) soit bijective,
pourvu que l'application
\begin{eqnarray*}
 \epsilon: H^0(\pp^n,\op{(n+1)}{n}(\ell))&\rightarrow&
\op{}{n-1}(\ell)_{Y_{y'+1}}\oplus\dots\oplus\op{}{n-1}(\ell)_{Y_{y'+f}}\oplus
\\
&& D\oplus \\
&&\tp{n}(\ell-1)_{Y_{y'+f+2}}\oplus\dots\oplus\tp{n}(\ell-1)_{Y_y}\oplus\\
&&\op{(n+1)}{n}(\ell)_{Z_{d+1}}\oplus\dots\oplus\op{(n+1)}{n}(\ell)_{Z_z}
\end{eqnarray*}
le soit. En utilisant alors un argument similaire {\`a} celui du
troisi{\`e}me cas de la preuve de (vi) on r{\'e}duit cette derni{\`e}re
assertion {\`a} la preuve de l'{\'e}nonc{\'e}
$$
{\bf RB}(\tp{n-1}(\ell-1),\op{n-1}(\ell),o_n(\ell-1)+y-y'-f-1,f;0,\delta)
$$

Dans les deux cas $g=0$ et $g\not=0$, le lecteur pourra v{\'e}rifier, pour que
ce qui pr{\'e}c{\`e}de ait un sens, que $y\geq y'+f+1$, ({\`a} noter que les
hypoth{\`e}ses $n\geq 2$ et $\ell\geq 0$ sont n{\'e}cessaires pour ce
faire). Ceci conclut donc la preuve du th{\'e}or{\`e}me~\ref{thmprinc}.\eop

\end{document}